\documentclass[11pt]{amsart}
\usepackage{amsmath, amsthm, amssymb, hyperref, geometry, url, mathrsfs}

\evensidemargin \oddsidemargin
\author{Benjamin Linowitz}
\address{Department of Mathematics\\University of Michigan\\Ann Arbor, MI 48109}
\email{linowitz@umich.edu}

\author{Jeffrey S. Meyer}
\address{Department of Mathematics\\
University of Oklahoma\\
Norman, OK 73019 USA}
\email[]{jmeyer@math.ou.edu}

\title{Systolic Surfaces of Arithmetic Hyperbolic 3-Manifolds}

\usepackage{fancyhdr}

\fancypagestyle{plain}

\DeclareMathAlphabet{\curly}{U}{rsfs}{m}{n}

\DeclareMathOperator{\Ann}{Ann}

\DeclareMathOperator{\covol}{covol}

\DeclareMathOperator{\disc}{disc}

\DeclareMathOperator{\Ram}{Ram}

\DeclareMathOperator{\PSL}{PSL}

\DeclareMathOperator{\sys}{sys}

\DeclareMathOperator{\Vol}{Vol}
\DeclareMathOperator{\coarea}{coarea}

\newtheorem{thm}{Theorem}[section]
\newtheorem{cor}[thm]{Corollary}

\theoremstyle{definition}

\newtheorem*{rmk}{Remark}

\theoremstyle{remark}

\newtheorem{proposition}{Proposition}[section]
\def\1{\mathbf{1}}
\def\disc{\mathrm{disc}}

\theoremstyle{plain}
\newtheorem{mainthm}{Theorem}

\theoremstyle{remark}

\theoremstyle{plain}
\newtheorem{theorem}[proposition]{Theorem}

\newtheorem{lemma}[proposition]{Lemma}

\def\1{\mathbf{1}}
\def\disc{\mathrm{disc}}

\newcommand{\abs}[1]{\left\vert#1\right\vert}

\newcommand{\iny}{\infty}

\makeatletter
\def\moverlay{\mathpalette\mov@rlay}
\def\mov@rlay#1#2{\leavevmode\vtop{%
   \baselineskip\z@skip \lineskiplimit-\maxdimen
   \ialign{\hfil$\m@th#1##$\hfil\cr#2\crcr}}}
\newcommand{\charfusion}[3][\mathord]{
    #1{\ifx#1\mathop\vphantom{#2}\fi
        \mathpalette\mov@rlay{#2\cr#3}
      }
    \ifx#1\mathop\expandafter\displaylimits\fi}
\makeatother

\makeatletter
\let\@@pmod\pmod
\DeclareRobustCommand{\pmod}{\@ifstar\@pmods\@@pmod}
\def\@pmods#1{\mkern4mu({\operator@font mod}\mkern 6mu#1)}
\makeatother

\begin{document}

\begin{abstract}
In this paper we examine the geometry of minimal surfaces of arithmetic hyperbolic $3$-manifolds.
In particular, we give bounds on the totally geodesic $2$-systole, construct infinitely many incommensurable manifolds with the same initial geometric genus spectrum in which volume and $1$-systole are controlled, and analyze the growth of the genera of minimal surfaces across commensurability classes.
These results have applications to the study of how Heegard genus grows across commensurability classes. 
\end{abstract}

\maketitle


\vspace{-2pc}

\section{Introduction}
Recent years have seen significant research in the geometry of $\pi_1$-injective surfaces of hyperbolic $3$-manifolds.  
The existence of $\pi_1$-injective surfaces was initially established for cusped finite volume hyperbolic $3$-manifolds in \cite{CooperLongReid}, for closed arithmetic hyperbolic 3-manifolds in \cite{L10}, and then finally for all closed hyperbolic $3$-manifolds in \cite{KahnMarkovic}.  
These results have had profound implications, in particular the resolution of Thurston's virtual Haken conjecture  \cite{Agol, BergeronWise}.
A natural continuation of this line of inquiry is to give an explicit description of the types of $\pi_1$-injective surfaces that can appear within a given hyperbolic $3$-manifold $M$.

The goal of this paper is to understand the minimal $\pi_1$-injective surfaces of $M$ when volume and $1$-systole are controlled.
We consider three distinct types of minimality.
The \textbf{systolic genus} of $M$, denoted $\mathrm{sysg}(M)$, is the minimal genus of a $\pi_1$-injective surface. 
The \textbf{$2$-systole} of $M$, denoted $\mathrm{sys}_2(M)$, is the minimal area of a $\pi_1$-injective surface.
The \textbf{totally geodesic $2$-systole} of $M$, denoted $\mathrm{sys}_2^{TG}(M)$, is the minimal area of an immersed totally geodesic surface.
Since hyperbolic $3$-manifolds are atoroidal, $\mathrm{sysg}(M)\ge2$. Furthermore, as totally geodesic surfaces are $\pi_1$-injective, $ \mathrm{sys}_2(M)\le\mathrm{sys}_2^{TG}(M)$. Finally, we denote by $\mathrm{sys}_1(M)$ the \textbf{$1$-systole} of $M$; that is, the least length of a closed geodesic on $M$.

Before stating our first result we introduce a piece of notation. 
Given two commensurable Kleinian groups $\Gamma_1$ and $\Gamma_2$, we define the \textbf{generalized index} of $\Gamma_2$ in $\Gamma_1$ to be $[\Gamma_1:\Gamma_2]:=\frac{[\Gamma_1:\Gamma_1\cap\Gamma_2]}{[\Gamma_2:\Gamma_1\cap\Gamma_2]}$. 
For a commensurability class of arithmetic Kleinian groups, we denote by $\Gamma_{\mathcal O}$ a group with minimal covolume in the class and by $\Gamma_{\mathcal O}^1$ a group contained in the class which is derived from the elements of reduced norm one in a maximal order contained in the class' invariant quaternion algebra. These groups and their construction will be described in detail in Section \ref{section:volumes}.

\begin{mainthm}\label{thm:volumebound}
Let $M$ be an arithmetic hyperbolic $3$-manifold which contains a finite area, properly immersed, totally geodesic surface.
Then $$\mathrm{sys}_2^{TG}(M) > c\left(\frac{\mathrm{Vol}(M)}{[\Gamma_{\mathcal O}:\pi_1(M)]}\right)^{1/2}$$ where $c>0$ is a constant depending only on the invariant trace field of $M$. If $\pi_1(M)$ contains a subgroup of the form $\Gamma_{\mathcal{O}}^1$ then we additionally have the inequality $$\mathrm{sys}_2^{TG}(M) < c' \mathrm{Vol}(M)^{1/2},$$ where $c'>0$ is a constant depending only on the invariant trace field of $M$.
\end{mainthm}

An immediate consequence of Theorem \ref{thm:volumebound} is that if $\pi_1(M)=\Gamma_{\mathcal O}$ then $\mathrm{sys}_2^{TG}(M)$ has order of magnitude $\mathrm{Vol}(M)^{1/2}$. In fact, because $[\Gamma_{\mathcal O}:\Gamma_{\mathcal O}^1]$ depends only on the groups' commensurability class, $\mathrm{sys}_2^{TG}(M)$ has order of magnitude $\mathrm{Vol}(M)^{1/2}$ whenever \[\Gamma_{\mathcal O} \supseteq \pi_1(M) \supseteq \Gamma_{\mathcal O}^1.\]

\begin{rmk}
If $\pi_1(M)$ is a maximal arithmetic group or a congruence subgroup of a maximal arithmetic group then the generalized index appearing in Theorem \ref{thm:volumebound} can easily be made explicit using Borel's classification of maximal arithmetic subgroups of $\PSL_2(\textbf{C})$ and the computation of their covolumes \cite{borel-commensurability}.
\end{rmk}

The proof of Theorem \ref{thm:volumebound} makes extensive use of Borel's computation of the covolumes of maximal arithmetic subgroups of $\PSL_2(\textbf{R})$ and $\PSL_2(\textbf{C})$ and the close relationship between the invariant quaternion algebra of an arithmetic hyperbolic $3$-manifold and those of its totally geodesic surfaces. An easy consequence of the proof of Theorem \ref{thm:volumebound} is the following.

\begin{cor}\label{cor:large}
For every $X>0$ there exist infinitely many arithmetic hyperbolic $3$-manifolds with $\mathrm{sys}_2^{TG}(M)>X$.
\end{cor}

Let $M$ be a closed hyperbolic $3$-manifold.
To each free homotopy class $[S]$ in $M$ of a $\pi_1$-injective, properly immersed, closed surface,
there is a well-defined minimal area, which we denote $\mathrm{area}(S)$ \cite{SchoenYau,SacksUhlenbeck}.
We define the \textbf{area set} of $M$ to be the set $\mathrm{A}(M)$ of all such areas.
For such a $\pi_1$-injective surface $S$ of genus $g$, the inequalities of Uhlenbeck and Thurston (cf., \cite[Lemma 6]{hass}) give
$$2\pi(g-1)\le \mathrm{area}(S)\le 4\pi(g-1).$$
By the work of Thurston, \cite[8.8.6]{Thurston} there are only finitely many $\pi_1$-injective, properly immersed, closed surfaces up to free homotopy of a fixed genus.
It follows that the area set is a discrete subset of the real numbers.
Furthermore, each area can only be realized by finitely many free homotopy classes of surfaces.

Every totally geodesic surface is the area minimizing representative in its homotopy class, and we denote the set of all areas contributed by totally geodesic surfaces $\mathrm{A}^{TG}(M)$. 
An interesting, though likely difficult, problem would be to study the distribution of $\mathrm{A}^{TG}(M)$ within $\mathrm{A}(M)$.
For $N\ge1$, we define $\mathrm{A}^{TG}(M)[N]$ to be the set of the first $N$ areas in $\mathrm{A}^{TG}(M)$.
The geometric genus set of $M$, denoted $\mathrm{GS}(M)$, is the set of isometry classes of finite area, properly immersed, totally geodesic surfaces of $M$ up to free homotopy. 
We define \[ \mathrm{GS}(M)[N]:=\{ Y\in\mathrm{GS}(M) : \mathrm{area}(Y)\in  \mathrm{A}^{TG}(M)[N] \}. \] 

Our next theorem revisits a construction of pairs of incommensurable arithmetic hyperbolic $3$-manifolds with the same set of commensurability classes of surfaces up to an arbitrary threshold that appeared in \cite[Section 7]{LMPT}. In particular we construct an infinite family of incommensurable arithmetic hyperbolic $3$-manifolds whose geometric genus spectra start off with the same $N$ terms and in which both volume and systole length are controlled.

\begin{mainthm}\label{thm:family}
Let $N\geq 1$. There exists an infinite sequence of incommensurable arithmetic hyperbolic $3$-manifolds $M_1,M_2,\dots$ such that $\mathrm{GS}(M_i)[N]=\mathrm{GS}(M_j)[N]\neq \emptyset$ for all $i,j$. In particular $\sys_2^{TG}(M_i)=\sys_2^{TG}(M_j)$ for all $i,j$. Furthermore, we have the estimates \[\Vol(M_n) < c_1\left(n\log(2n)\right)^{3/2}, \qquad \mbox{ and } \qquad \sys_1(M_n) < c_2\log(n),\] where $c_1,c_2$ are positive constants depending on $N$.
\end{mainthm}

Theorem \ref{thm:family} may be viewed as a higher dimensional analogue of a recent result of Millichap \cite{Mill}, who exhibited an infinite family of incommensurable hyperbolic $3$-manifolds, all of whose complex length spectra begin with the same $2n+1$ lengths.

It is natural to ask about the minimum value assumed by $\mathrm{sysg}(M)$, $\sys_2(M)$ or $\sys_2^{TG}(M)$ as $M$ varies over the set of arithmetic hyperbolic $3$-manifolds lying in a fixed commensurability class. Our next result will show that given any real number $x>0$, almost every commensurability class of arithmetic hyperbolic $3$-manifolds has the property that every representative $M$ satisfies $\mathrm{sysg}(M), \sys_2(M),\sys_2^{TG}(M)>x$. In order to make this precise we introduce some additional notation.

Let $k$ be a number field with a unique complex place and $\mathscr{C}$ be a commensurability class of arithmetic hyperbolic $3$-manifolds having invariant trace field $k$. In Section \ref{section:volumes} we will define a notion of volume on $\mathscr C$. This volume is in terms of the volume of a certain distinguished representative of the class which arises in a natural way in the construction of arithmetic lattices in $\PSL_2(\textbf{C})$. Let $N_k(V)$ denote the number of commensurability classes of arithmetic hyperbolic $3$-manifolds with volume less than $V$ and invariant trace field $k$. Additionally, let $N_k^g(V;x)$ (respectively $N_k^2(V;x)$ and $N_k^{TG}(V;x)$) denote the number of commensurability classes of arithmetic hyperbolic $3$-manifolds with invariant trace field $k$, volume less than $V$ and which contain a representative $M$ satisfying $\mathrm{sysg}(M)<x$ (respectively $\sys_2(M)<x$ and $\sys_2^{TG}(M)<x$).

\begin{mainthm}\label{thm:limit}
For all sufficiently large $x$ we have \[ \lim_{V \to \iny} \frac{N_k^g(V;x)}{N_k(V)} = \lim_{V \to \iny} \frac{N_k^2(V;x)}{N_k(V)} = \lim_{V \to \iny} \frac{N_k^{TG}(V;x)}{N_k(V)} = 0. \]
\end{mainthm}

The Heegard genus of a closed $3$-manifold, denoted $\mathrm{Hg}(M)$, is the minimal genus of a Heegard surface.
There has been considerable work towards understanding how Heegard genus increases along covers \cite{Belolipetsky,L06,GR09,DB}. 
In particular, in the case of congruence covers of a compact arithmetic 3-manifold,  \cite{L06} and \cite{GR09} independently have shown that the Heegard genus grows linearly with the degree of the covers.
Since the Heegard genus of a $3$-manifold is bounded below by its systolic genus, Theorem \ref{thm:limit} implies that the Heegard genus grows across commensurability classes.  
More precisely, let $\mathrm{Hg}_k(V;x)$ denote the number of commensurability classes of arithmetic hyperbolic $3$-manifolds with invariant trace field $k$, volume less than $V$ and which contain a representative $M$ satisfying $\mathrm{Hg}(M)<x$.

\begin{cor}\label{cor:Heegard}
For all sufficiently large $x$ we have  $\displaystyle\lim_{V \to \iny} \dfrac{\mathrm{Hg}_k(V;x)}{N_k(V)}= 0.$
\end{cor}


\section{Notation}

Let $k$ be a number field of degree $n_k$ and discriminant $\Delta_k$. We denote the ring of integers of $k$ by $\mathcal O_k$ and the Dedekind zeta function of $k$ by $\zeta_k(s)$. Given an ideal $I\subset \mathcal O_k$, we denote by $\abs{I}$ the norm of $I$.

If $A$ is a quaternion algebra over $k$ then we denote by $\Ram(A)$ (respectively, $\Ram_f(A)$) the set of primes (respectively, finite primes) of $k$ which ramify in $A$. The discriminant of $A$, denoted $\disc(A)$, is defined to be the product of all primes in $\Ram_f(A)$. 

We denote by $\textbf{H}^2$ and $\textbf{H}^3$ real hyperbolic spaces of dimensions $2$ and $3$, and will refer to lattices in $\PSL_2(\textbf{R})$ as Fuchsian and lattices in $\PSL_2(\textbf{C})$ as Kleinian.

\section{Arithmetic Fuchsian and Kleinian groups and their covolumes}\label{section:volumes}

In this section we briefly describe the construction of arithmetic lattices in $\PSL_2(\textbf{R})$ and $\PSL_2(\textbf{C})$. For a more detailed treatment we refer the reader to \cite{MR}.

Let $\ell$ be a totally real number field and $B$ a quaternion algebra over $\ell$ which is unramified at a unique real prime $v$. We therefore have an isomorphism $B\otimes_\ell \ell_v \cong \mathrm{M}_2(\textbf{R})$. Let $\mathcal O$ be a maximal order of $B$ and $\mathcal O^1$ be the subgroup of $\mathcal O^*$ consisting of those elements with reduced norm one. We denote by $\Gamma_{\mathcal O}^1$ the image in $\PSL_2(\textbf{R})$ of $\mathcal O^1$ and by $\Gamma_{\mathcal O}$ the image in $\PSL_2(\textbf{R})$ of $N(\mathcal O)$, where $N(\mathcal O)$ is the normalizer in $B^*$ of $\mathcal O$. The groups $\Gamma_{\mathcal O}$ and $\Gamma_{\mathcal O}^1$ are discrete subgroups of $\PSL_2(\textbf{R})$ which are cocompact whenever $B$ is a division algebra. The coareas of these groups are given by the following formulas (see Borel \cite{borel-commensurability} and Chinburg and Friedman  \cite[Proposition 2.1]{chinburg-smallestorbifold}) :

\begin{equation}\label{equation:fuchsiannorm1volume}
\coarea(\Gamma_\mathcal{O}^1) = \frac{8\pi\abs{\Delta_\ell}^{\frac{3}{2}}\zeta_\ell(2)}{(4\pi^2)^{n_\ell}}\prod_{\mathfrak p\in\Ram_f(B)}\left(\abs{\mathfrak p}-1\right), 
\end{equation}
and
\begin{equation}\label{equation:fuchsianvolume}
\coarea(\Gamma_\mathcal{O})= \frac{4\pi\abs{\Delta_\ell}^{\frac{3}{2}}\zeta_\ell(2)}{(4\pi^2)^{n_\ell}t(B)}\prod_{\mathfrak p\in\Ram_f(B)}\frac{\abs{\mathfrak p}-1}{2},
\end{equation}
where $t(B)$ is the type number of $B$. It is known that the group $\Gamma_{\mathcal O}$ has minimal covolume within its commensurability class. Any group which is commensurable with a group $\Gamma_{\mathcal O}^1$ as constructed above is called an \textbf{arithmetic Fuchsian group}.

We now describe the construction of arithmetic lattices in $\PSL_2(\textbf{C})$. Let $k$ be a number field with a unique complex place $v$ and $A$ a quaternion algebra over $k$ which is ramified at all real places of $k$. Let $\mathcal O$ be a maximal order of $A$ and define subgroups $\Gamma_\mathcal{O}^1$ and $\Gamma_\mathcal{O}$ of $\PSL_2(\textbf{C})$ (constructed as above) via the identification $A\otimes_k k_v\cong \mathrm{M}_2(\textbf{C})$. The groups $\Gamma_\mathcal{O}^1$ and $\Gamma_\mathcal{O}$ are discrete subgroups of $\PSL_2(\textbf{C})$ of finite covolume which are cocompact whenever $A\not\cong \mathrm{M}_2(\textbf{Q}(\sqrt{-d}))$ for some imaginary quadratic field $\textbf{Q}(\sqrt{-d})$. Furthermore, $\Gamma_\mathcal{O}$ has minimal covolume within its commensurability class. Any discrete subgroup of $\PSL_2(\textbf{C})$ which is commensurable with a group of the form $\Gamma_{\mathcal O}^1$ as constructed above is called an \textbf{arithmetic Kleinian group}. The work of Borel \cite{borel-commensurability} (see also \cite[Proposition 2.1]{chinburg-smallestorbifold}) shows that the covolumes of $\Gamma_\mathcal{O}^1$ and $\Gamma_\mathcal{O}$ are given by

\begin{equation}\label{equation:kleiniannorm1volume}
\coarea(\Gamma_\mathcal{O}^1) = \frac{\abs{\Delta_k}^{\frac{3}{2}}\zeta_k(2)}{(4\pi^2)^{n_k-1}}\prod_{\mathfrak p\in\Ram_f(A)}\left(\abs{\mathfrak p}-1\right),
\end{equation}
and 
\begin{equation}\label{equation:kleinianvolume}
\covol(\Gamma_{\mathcal O})= \frac{2\pi^2\abs{\Delta_k}^{\frac{3}{2}}\zeta_k(2)}{(4\pi^2)^{n_k}t(A)}\prod_{\mathfrak p\in\Ram_f(A)}\frac{\abs{\mathfrak p}-1}{2}.
\end{equation}

A pair of arithmetic lattices $\Gamma_1,\Gamma_2$ derived from quaternion algebras $(k_1,A_1), (k_2,A_2)$ will be commensurable (in the wide sense) precisely when $k_1\cong k_2$ and $A_1\cong A_2$ \cite[Theorem 8.4.1]{MR}. We will make use of this fact throughout this paper.

Given a commensurability class $\mathscr C$ of arithmetic hyperbolic $3$-manifolds arising from the arithmetic data $(k,A)$, we define the \textbf{volume} $V_\mathscr{C}$ of $\mathscr{C}$ to be $V_\mathscr{C}:=\covol(\Gamma_{\mathcal O}^1)$, where $\mathcal O$ is a maximal order of $A$. It follows from the work of Borel \cite{borel-commensurability} that there are only finite many commensurability classes of bounded volume.


\section{Proof of Theorem \ref{thm:volumebound}}

We now prove Theorem \ref{thm:volumebound}.

Let $\Gamma=\pi_1(M)$, and denote by $k=k\Gamma$ be the invariant trace field of $M$ and by $A=A\Gamma$ the invariant quaternion algebra of $M$. The results of \cite[Chapter 9.5]{MR} show that all totally geodesic surfaces in $M$ arise from quaternion algebras $B$ over $\ell$, where $\ell$ is the maximal totally real subfield of $k$ (and satisfies $[k:\ell]=2$) and $B$ is a quaternion algebra defined over $\ell$ which is unramified at a unique real place of $\ell$ (the one lying under the complex place of $k$) and has the property that $A\cong B\otimes_\ell k$. The requirement that $A\cong B\otimes_\ell k$ imposes severe restrictions on the ramification of $B$, as the following result \cite[Theorem 9.5.5]{MR} makes clear.

\begin{theorem}\label{theorem:tgs}
Let $k$ be a number field with a unique complex place and $A$ be a quaternion algebra over $k$ which is ramified at all real places. Let $k$ be such that $[k:\ell]=2$, where $\ell=k\cap \textbf{R}$. Let $B$ be a quaternion algebra over $\ell$ which is ramified at all real places except the identity. Then $A\cong B\otimes_\ell k$ if and only if $\Ram_f(A)$ consists of $2r$ distinct ideals $\{ \mathfrak P_1,\mathfrak P_1',\dots, \mathfrak P_r,\mathfrak P_r'\}$, where \[\mathfrak P_i\cap \mathcal O_\ell = \mathfrak P_i'\cap \mathcal O_\ell = \mathfrak p_i\] and $\Ram_f(B)\supseteq \{\mathfrak p_1,\dots,\mathfrak p_r\}$ with $\Ram_f(B)\setminus  \{\mathfrak p_1,\dots,\mathfrak p_r\}$ consists of primes of $\mathcal O_\ell$ which are either inert or ramified in the extension $k/\ell$.
\end{theorem}

As $\mathrm{GS}(M)\neq \emptyset$, the set $\Ram_f(A)$ is of the form $\{ \mathfrak P_1,\mathfrak P_1',\dots, \mathfrak P_r,\mathfrak P_r'\}$ described in Theorem \ref{theorem:tgs}. Every non-elementary Fuchsian subgroup of $\Gamma$ is contained in an arithmetic Fuchsian group \cite[Theorem 9.5.2]{MR} which arises from a quaternion algebra $B$ over $\ell$ which is ramified at $\{\mathfrak p_1,\dots,\mathfrak p_r\}$ by Theorem \ref{theorem:tgs}. Note that if $B$ is an indefinite quaternion algebra over $\ell$ then the type number $t(B)$ of $B$ is equal to $[\ell(B):\ell]$ where $\ell(B)$ is the maximal abelian extension of $\ell$ of exponent $2$ which is unramified outside the real places in $\Ram(B)$ and in which all primes of $\Ram_f(B)$ split completely (see \cite[p. 39]{Chinburg-Friedman-selectivity}). Observe that $\ell(B)$ is contained in the narrow class field of $\ell$, which has degree over $\ell$ at most $2^{n_\ell}h_\ell$. Section \ref{section:volumes} and equation (\ref{equation:fuchsianvolume}) now show that 
\begin{equation}\label{equation:sysbound1}
\mathrm{sys}_2^{TG}(M) \geq \frac{4\pi\abs{\Delta_\ell}^{\frac{3}{2}}\zeta_\ell(2)}{(16\pi^2)^{n_\ell}h_\ell}\prod_{i=1}^r \frac{\abs{\mathfrak p_i}-1}{2}=c_{\ell}\prod_{i=1}^r \frac{\abs{\mathfrak p_i}-1}{2},
\end{equation}
where $c_\ell>0$ is a constant depending only on $\ell$ (and hence on $k$ as $\ell=k\cap \textbf{R}$ is the maximal totally real subfield of $k$).

Recall that $ \Ram_f(A) = \{ \mathfrak P_1,\mathfrak P_1',\dots, \mathfrak P_r,\mathfrak P_r'\}$ and $\abs{\mathfrak P_i}=\abs{\mathfrak P_i'}=\abs{\mathfrak p_i}$ for $i=1,\dots,r$. It follows that \[\prod_{\substack{\mathfrak p\in\Ram_f(A)\\ \abs{\mathfrak p}>2}}\frac{\abs{\mathfrak p}-1}{2} \geq \left(\prod_{i=1}^r \frac{\abs{\mathfrak p_i}-1}{2}\right)^2,\] hence the first part of the theorem follows from equations (\ref{equation:sysbound1}) and (\ref{equation:kleinianvolume}).

We now suppose that $\pi_1(M)$ contains a subgroup of the form $\Gamma_{\mathcal O}^1$ and will derive an upper bound for $\mathrm{sys}_2^{TG}(M)$. That $\mathrm{GS}(M)\neq \emptyset$ implies that $ \Ram_f(A) = \{ \mathfrak P_1,\mathfrak P_1',\dots, \mathfrak P_r,\mathfrak P_r'\}$ where the primes $\mathfrak P_i$ and $\mathfrak P_i'$ both lie above the prime $\mathfrak p_i$ of $\mathcal O_\ell$. Let $\mathfrak p$ be a prime of $\mathcal O_\ell$ which is inert in the extension $k/\ell$. The effective Chebotarev density theorem \cite{effectiveCDT} shows that there exists an absolute constant $C>0$ such that we may assume $\abs{\mathfrak p}< \abs{\Delta_\ell}^C$. We now define a quaternion algebra $B$ over $\ell$ by declaring it to be ramified at all real places except the identity and at the finite primes $\{\mathfrak p_1,\dots,\mathfrak p_r\}$ if $r+n_\ell-1$ is even and $\{ \mathfrak p_1,\dots,\mathfrak p_r\}\cup\{\mathfrak p\}$ if $r+n_\ell-1$ is odd. Theorem \ref{theorem:tgs} shows that $A\cong B\otimes_\ell k$. Let $\mathcal O_B$ be a maximal order of $B$, chosen so that $\mathcal O_B\otimes_{\mathcal O_\ell}\mathcal O_k\cong \mathcal O$. Therefore $\Gamma_{\mathcal O_B}^1\subset \Gamma_{\mathcal O}^1\subset \pi_1(M)$ and $\mathrm{sys}_2^{TG}(M)\leq \coarea(\Gamma_{\mathcal O_B}^1)$. Borel \cite{borel-commensurability} (see also \cite[Chapter 11.6]{MR}) has shown that $[\Gamma_{\mathcal O}:\Gamma_\mathcal O^1]$ depends only on $k$, hence it suffices to show that $\coarea(\Gamma_{\mathcal O_B}^1)\leq c_k \covol(\Gamma_{\mathcal O}^1)^{\frac{1}{2}}$ for some constant $c_k>0$ depending only on $k$. Notice that $\abs{\mathfrak p}$ has been bounded by a quantity which depends only on $\ell$, and that \[\left(\prod_{i=1}^r(\abs{{\mathfrak p_i}}-1)\right)^2=\prod_{i=1}^r(\abs{\mathfrak P_i}-1)(\abs{\mathfrak P_i'}-1)=\prod_{\mathfrak q\in\Ram_f(A)}(\abs{\mathfrak q}-1).\] The upper bound for $\mathrm{sys}_2^{TG}(M)$ now follows from equations (\ref{equation:fuchsiannorm1volume}) and (\ref{equation:kleiniannorm1volume}).

We conclude this section by proving Corollary \ref{cor:large}. 

Let $k$ be a number field with a unique complex place and whose maximal totally real subfield $\ell$ satisfies $[k:\ell]=2$. Infinitely many primes of $\ell$ split in the quadratic extension $k/\ell$, hence there exist finitely many primes $\mathfrak p$ of $\mathcal O_\ell$ such that \[ \abs{\mathfrak p} > X\left(\frac{(8\pi^2)^{n_\ell}h_\ell\abs{\Delta_\ell}^{-\frac{3}{2}}}{2\pi\zeta_\ell(2)}\right)+1.\] Fix one such prime $\mathfrak p$. Let $A$ be the quaternion algebra over $k$ which is ramified at all real primes of $k$ and at the two finite primes lying above $\mathfrak p$. Let $\Gamma$ be any torsion-free arithmetic Kleinian group arising from $(A,k)$. The results of \cite[Chapter 9.5]{MR} show that any non-elementary Fuchsian subgroup of $\Gamma$ must arise from an indefinite quaternion algebra $B$ over $\ell$ which is ramified at $\mathfrak p$ (and potentially other primes). Note that if $t(B)$ is the type number of an indefinite quaternion algebra over $\ell$ then we have already seen that $t(B)\leq 2^{n_\ell}h_\ell$. The corollary now follows from Section \ref{section:volumes} and in particular equation (\ref{equation:fuchsianvolume}).


\section{Proof of Theorem \ref{thm:family}}

Our proof will make use of the following refinement of Linnik's theorem on primes in arithmetic progressions, the proof of which follows from \cite[Corollary 18.8, p. 442]{IK} (see also \cite[Proposition 8.1]{LMP}).

\begin{lemma}\label{lemma:linnik}
There is an absolute constant $L>0$ such that for every pair of coprime integers $a,q$ with $q\geq 2$ the $n$-th prime $p\equiv a\pmod{q}$ satisfies \[p\leq q^L\cdot n\log(2n).\]
\end{lemma}

Let $m$ be a natural number such that the $m$-th prime $p_m$ satisfies $p_m>2^{10} 3^2 p_{N+1}$.

Let $L_1, L_2, \dots$ be a sequence of imaginary quadratic fields of discriminants $-q_1, -q_2,\dots$ (where all of the $q_i$ are prime) in which the rational prime $13$ splits and in which $p_1,p_2,p_3,p_4,p_5,p_7,\dots,p_{m-1}$ are inert. We assume as well that none of the $L_i$ are contained in one of the cyclotomic fields $\textbf{Q}(\sqrt{-1}),\textbf{Q}(\sqrt{-3}), \textbf{Q}(\zeta_5)$. Note that these fields may be constructed by selecting the primes $q_1,q_2,\dots$ to lie within a certain residue class modulo $C=4p_1\cdots p_{m-1}$. It now follows from Lemma \ref{lemma:linnik} that we may take \begin{equation}\label{equation:discbound} \abs{\Delta_{L_i}} < C^L \cdot i\log(2i)\end{equation} for some absolute constant $L>0$.

For every $i\geq 1$ let $A_i$ be the quaternion algebra over $L_i$ such that $\disc(A_i)=13\mathcal O_{L_i}$. Let $\mathcal O_i$ be a maximal order of $A_i$ and observe that $\Gamma_{\mathcal O_i}$ is torsion-free as no quadratic cyclotomic extension of $L_i$ embeds into $A_i$ (note that $13$ splits in both $\textbf{Q}(i)/\textbf{Q}$ and $\textbf{Q}(\sqrt{-3})/\textbf{Q}$, hence neither $L_j(i)$ nor $L_j(\sqrt{-3})$ embed into $A_j$).

A consequence of Theorem \ref{theorem:tgs} is that if a quaternion algebra $B$ over $\textbf{Q}$ satisfies $B\otimes_{\textbf{Q}}L_i\cong A_i$ but $B\otimes_{\textbf{Q}}L_j\not\cong A_j$ is that $\Ram(B)$ must contain a prime $p$ which splits in $L_i$ but not $L_j$. By construction this prime must satisfy $p\geq p_m$. Let $B_1, B_2, \dots, B_5, B_7,\dots, B_{N+1}$ be indefinite rational quaternion algebras for which $\Ram(B_i)=\{ 13,p_i \}$ and $\mathcal O_{B_i}$ be a maximal order of $B_i$. We now deduce from equations (\ref{equation:fuchsiannorm1volume}) and (\ref{equation:fuchsianvolume}) that every arithmetic Fuchsian derived from $B$ has coarea larger than all of the $\Gamma_{\mathcal O_{B_i}}$. Because the $B_i$ embed into all of the quaternion algebras $A_1, A_2, \dots$, we see that $\mathrm{GS}(\textbf{H}^3/\Gamma_{\mathcal O_i})[N]=\mathrm{GS}(\textbf{H}^3/\Gamma_{\mathcal O_j})[N]$ for all $i,j$. The assertion about $\Vol(M_n)$ now follows from equations (\ref{equation:discbound}) and (\ref{equation:kleiniannorm1volume}), along with the fact that $\zeta_{L_i}(2)<\zeta(2)^2<3$ for all $i$.

To conclude the proof of Theorem \ref{thm:family}, we must derive an upper bound for $\sys_1(\textbf{H}^3/\Gamma_{\mathcal O_i})$. This upper bound is an immediate consequence of the well-known fact that if $M$ is a closed hyperbolic $3$-manifold then $\sys_1(M)<\log(\Vol(M))$ (this follows from the injectivity radius of $M$ being $\frac{1}{2}\sys_1(M)$ and the fact that the volume of a ball of radius $r$ in $\textbf{H}^3$ has order of magnitude $e^r$).


\section{Proof of Theorem \ref{thm:limit}}	
	
If $M$ is a closed hyperbolic $3$-manifold then it is clear that $\sys_2^{TG}(M)\geq \sys_2(M)$, hence it suffices to prove Theorem \ref{thm:limit} for $N_k^g(V;x)$ and $N_k^2(V;x)$.
	
We begin by proving the theorem for $N_k^g(V;x)$. In essence, our proof follows from two previously known results. The first is due to the first author and McReynolds, Pollack and Thompson \cite{LMPT2} and concerns the number of commensurability classes of arithmetic hyperbolic $3$-manifolds containing a representative with systole less than any fixed bound. The second result is an important systole inequality of Belolipetsky \cite{Belolipetsky}, whose proof relies on a high genus systole inequality of Gromov \cite{GR83} (see also \cite[p. 88]{K07}):
 	
\begin{theorem}[Belolipetsky] Let $M$ be a closed hyperbolic $3$-manifold. For any $\epsilon>0$, if the systole $\sys_1(M)$ of $M$ is sufficiently large then \begin{equation}\label{equation:misha}\mathrm{sysg}(M)\geq e^{(\frac{1}{2}-\epsilon)\sys_1(M)}.\end{equation}
\end{theorem}

Let $x_0=x_0(\epsilon)$ be large enough so that (\ref{equation:misha}) holds with $\epsilon=\frac{1}{4}$ whenever $\sys_1(M)>x_0$. If $x>x_0$ and $\mathrm{sysg}(M)<x$ then (\ref{equation:misha}) shows that $\sys_1(M)<4\log(x)$. Given a function $f(x)$, define $N_k(V;f(x))$ to be the number of commensurability classes of arithmetic hyperbolic $3$-manifolds with invariant trace field $k$, volume less than $V$ and which contain a representative $M$ with systole $\sys_1(M)<f(x)$. We now see that \[N_k^g(V;x)< N_k(V;x_0)+N_k(V; 4\log(x)),\] as for a given commensurability class $\mathscr C$, either $\mathscr C$ contains a manifold $M$ with $\sys_1(M)\leq x_0$ in which case (\ref{equation:misha}) need not apply, or else every manifold $M\in\mathscr C$ satisfies \[ \mathrm{sysg}(M)<x \Longrightarrow x_0 < \sys_1(M) \leq 4\log(x). \] Our proof in the case of $N_k^g(V;x)$ now follows from \cite[Corollary 1.2]{LMPT2}. The remaining case of $N_k^2(V;x)$ follows from exactly the same argument, replacing the aforementioned systole inequality of Belolipetsky by \cite[Corollary 2.3]{Belolipetsky}.


\subsection*{Acknowledgments} The first author was partially supported by an NSF RTG grant DMS-1045119 and an NSF Mathematical Sciences Postdoctoral Fellowship.



\end{document}